\documentclass{article}
\usepackage{amsmath,amsgen,amstext,amsbsy,amsopn,amsfonts,graphics,graphicx}
\begin{document}

\centerline{\LARGE Triacontagonal coordinates for the $E_8$ root system}

\begin{flushleft}
David A. Richter \\
Department of Mathematics, MS 5248 \\
Western Michigan University \\
Kalamazoo MI 49008-5248 \\
\texttt{david.richter@wmich.edu} \\
\end{flushleft}

Abstract. This note gives an explicit formula for the elements of the $E_8$ root system.  The formula is triacontagonally
symmetric in that one may clearly see an action by the cyclic group with 30 elements.  The existence of such a
formula is due to the fact that the Coxeter number of $E_8$ is 30.



\section{Introduction}

The Coxeter numbers of $H_4$ and $E_8$ are both equal to 30, \cite{Hum}.  Some artistically-inclined mathematicians
have used this fact in order to depict these root systems.  In so doing, each has necessarily produced
a figure having triacontagonal symmetry, meaning the same symmetry as a regular 30-sided polygon.  
According to \cite{Cox1}, where it is used as the frontispiece, 
van Oss first sketched such a projection of the regular polytope having Schl\"afli symbol $\{3,3,5\}$, also known
as the ``600-cell''.  The vertices of $\{3,3,5\}$
coincide with the elements of the (non-crystallographic) root system $H_4$.  One may also find this sketch in
the article \cite{Sti} and a sketch of the dual polytope $\{5,3,3\}$ in \cite{Chi}.

The roots of $E_8$ coincide with the vertices of a highly symmetric convex polytope $4_{21}$ apparently discovered by 
Gosset, \cite{Cox1,Gos}, and also with the vertices of a regular complex polytope $3\{3\}3\{3\}3\{3\}3$
discovered by Witting, \cite{Cox2,CoSh}.  The book \cite{Cox2} displays a triacontagonal projection of 
$3\{3\}3\{3\}3\{3\}3$ and \cite{CoSh} displays a triacontagonal projection of $4_{21}$.
According to \cite{McM}, Peter McMullen sketched the triacontagonal projections of both of these polytopes
by hand during a period lasting from August 18 until August 20 in 1964.
Recently John Stembridge produced an 8-color computer sketch of $4_{21}$, to supplement an
announcement by a research group at the American Institute of Mathematics of having computed a particular
class of representations of a real Lie group having the $E_8$ root system.  A simplified sketch is offered here,
showing only the relative locations of the 240 vectors, in 8 concentric cycles with 30 points per cycle.

Here is a streamlined version of the formula from \cite{Cox1} for the elements of the $H_4$ root system.
Let $\omega=\exp\left(\frac{i\pi}{30}\right)$ and $a>b>c>d$ be the positive
roots of the polynomial $45x^8-90x^6+60x^4-15x^2+1$ over $\mathbb{R}$.  Then the 120 vectors
$$\begin{array}{ll}
A_n = (a\omega^{2n},d\omega^{22n}), & B_n = (b\omega^{2n+1},c\omega^{22n+11}), \\
C_n = (c\omega^{2n+1},-b\omega^{22n+11}), & D_n = (d\omega^{2n},-a\omega^{22n}), \\
\end{array}$$
where $n\in\{0,1,2,...,28,29\}$, comprise the $H_4$ root system.
This formula has the feature that projecting along the first coordinate,
essentially forgetting the second coordinate,
yields the triacontagonal projection of the $H_4$ roots.
The purpose of this note is to give a similar formula for the $E_8$ root system.

The triacontagonal projections of the $H_4$ and $E_8$ root systems are closely related.  
The latter is the union of two scaled sizes of the former, where the
ratio of the larger to the smaller is the golden ratio $\tau=\frac{1}{2}\left(1+\sqrt{5}\right)$.
Given this, one expects that a similar formula should exist for the $E_8$ roots.
Indeed, the formula given here was obtained by ``guess-and-check''.  Checking is tedious, but
one may use Maple, for example, to compute as many inner products as are necessary
to become convinced that the formula yields a root system isomorphic to $E_8$.
(A more 	``ethical'' way to obtain the formula would be to diagonalize a Coxeter element
of the Coxeter group of $E_8$.  However, this author believes that the ends could not possibly justify
the enormity of the computations involved in such means.)

\begin{centering}
\includegraphics[width=5 in]{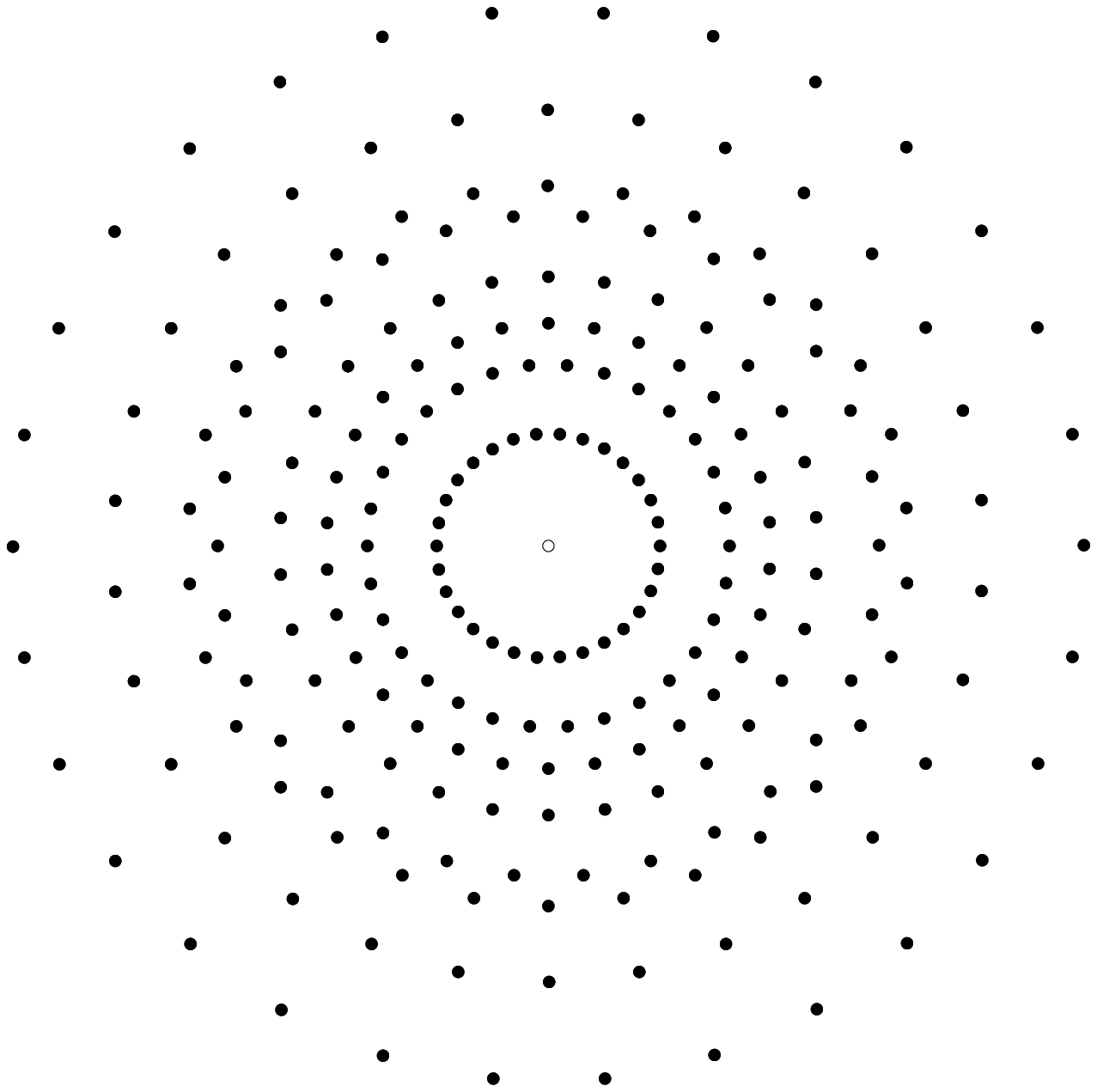} \\
{\bf Figure 1.}  Triacontagonal projection of the $E_8$ root system. \\
\end{centering}

\section{The formula}


Denote $\omega=\exp\left(\frac{i\pi}{30}\right)$ and define $\{a,b,c,d\}$ as above.
More explicitly, $a$, $b$, $c$, and $d$ are the positive numbers satisfying
$$\begin{array}{ll}
2a^2=1+3^{-1/2}5^{-1/4}\tau^{3/2}, & 2b^2=1+3^{-1/2}5^{-1/4}\tau^{-3/2}, \\
2c^2=1-3^{-1/2}5^{-1/4}\tau^{-3/2}, & 2d^2=1-3^{-1/2}5^{-1/4}\tau^{3/2}, \\
\end{array}$$
where $\tau$ is the golden ratio.  For any integer $n$, denote
$c_n=\omega^n+\omega^{-n}=2\cos\left(\frac{n\pi}{30}\right)$.  Next, denote
$$\begin{array}{llll}
r_1 = a/c_9, &
r_2 = b/c_9, &
r_3 = c/c_9, &
r_4 = d/c_9, \\
r_5 = a/c_3, &
r_6 = b/c_3, &
r_7 = c/c_3, &
r_8 = d/c_3. \\
\end{array}$$
Then the 240 rows of the matrices
$$\left[\begin{array}{c}
A_n \\
B_n \\
C_n \\
D_n \\
E_n \\
F_n \\
G_n \\
H_n \\
\end{array}\right] =
\left[\begin{array}{cccc}
 r_1         &  r_4         &  r_6\omega      &  r_7\omega  \\
 r_2\omega^{29} &  r_3\omega^{19} & -r_8\omega^{24} & -r_5\omega^{18} \\
 r_3\omega^{29} & -r_2\omega^{19} &  r_5\omega^{24} & -r_8\omega^{18} \\
 r_4         & -r_1         &  r_7\omega      & -r_6\omega  \\
 r_5         &  r_8         & -r_2\omega      & -r_3\omega  \\
 r_6\omega^{29} &  r_7\omega^{19} &  r_4\omega^{24} &  r_1\omega^{18} \\
 r_7\omega^{29} & -r_6\omega^{19} & -r_1\omega^{24} &  r_4\omega^{18} \\
 r_8         & -r_5         & -r_3\omega      &  r_2\omega  \\
\end{array}\right]\cdot
\left[\begin{array}{cccc}
\omega^{2n} & & & \\
 & \omega^{22n} & & \\
 & & \omega^{14n} & \\
 & & & \omega^{26n} \\
\end{array}\right],$$
where $n\in\{0,1,2,...,28,29\}$, as regarded as elements of $\mathbb{C}^4$,
comprise a root system isomorphic to $E_8$.  Moreover, each of these
vectors has norm equal to 1.


\section{Remarks}

Projecting along the first coordinate, by forgetting the last three coordinates, yields the image depicted above as a
subset of $\mathbb{C}\cong\mathbb{R}^2$.

Each cycle of 30 roots, as denoted by $A_n$, $B_n$, and so on, is expressed using the trigonometric function
$n\mapsto(\omega^{2n},\omega^{22n},\omega^{14n},\omega^{26n})$, composed with some amplitude and phase adjustments. 
Each phase shift was chosen so that $\{A_n,B_n,C_n,D_n,E_n,F_n,G_n,H_n\}$ is a system
of simple roots for any fixed value of $n$.

\begin{centering}
\begin{picture}(200,100)
\put(25,25){\circle*{5}}
\put(75,25){\circle*{5}}
\put(125,25){\circle*{5}}
\put(175,25){\circle*{5}}
\put(25,75){\circle*{5}}
\put(75,75){\circle*{5}}
\put(125,75){\circle*{5}}
\put(175,75){\circle*{5}}
\put(25,25){\line(1,0){150}}
\put(25,75){\line(1,0){150}}
\put(75,25){\line(-1,1){50}}
\put(20,12){$F_n$}
\put(70,12){$A_n$}
\put(120,12){$C_n$}
\put(170,12){$D_n$}
\put(20,81){$B_n$}
\put(70,81){$E_n$}
\put(120,81){$G_n$}
\put(170,81){$H_n$}
\end{picture} \\
{\bf Figure 2.}  The Dynkin diagram of $E_8$. \\
\end{centering}
\vskip .25 in

The amplitudes $r_k$ were chosen so that each vector has norm 1.
Alternatively, one may use the amplitudes
$$\begin{array}{llll}
r_1 = 1, & 
r_2 = c_{11}, & 
r_3 = c_6c_{13}, & 
r_4 = c_6c_{14}, \\
r_5 = c_{12}, & 
r_6 = c_{11}c_{12}, & 
r_7 = c_{13}, & 
r_8 = c_{14}. \\
\end{array}$$
In so doing, one avoids the complicated definition/formulae for the values $\{a,b,c,d\}$ given above.
For example, using this choice of amplitudes facilitates the verification of the formula, for then all the coordinates
lie in the cyclotomic field $\mathbb{Q}(\omega)$.
The tradeoff is that the norms of the vectors are more complicated to express.


\end{document}